\documentclass[11pt]{article}
\usepackage{amsmath,amstext,amssymb}
\usepackage{setspace,graphicx}
\usepackage{wrapfig}
\title{A Generalization \\of the B\^{o}cher-Grace Theorem}
\author{John Clifford and Michael Lachance}

\newcommand{\sech}{\mbox{\text{sech }}}

 \begin{spacing}{1}
\begin{document}

\maketitle \abstract\noindent{The B\^{o}cher-Grace
Theorem can be stated as follows:  Let $p$ be a third degree complex
polynomial. Then there is a unique inscribed ellipse interpolating
the midpoints of the triangle formed from the roots of $p$, and the
foci of the ellipse are the critical points of $p$.  Here, we prove
the following generalization:  Let $p$ be an $n^{th}$ degree complex
polynomial and let its critical points take the form
$$
\alpha+\beta \cos k\pi/n, \quad k=1,\ldots,n-1, \quad\beta\ne0.
$$
Then there is an inscribed ellipse interpolating the midpoints of
the convex polygon formed by the roots of $p$, and the foci of this
ellipse are the two most extreme critical points of $p$:
~$\alpha\pm\beta \cos \pi/n$. }
\section{Introduction}

The B\^{o}cher-Grace Theorem has been discovered
independently by many mathematicians.  Recently, proofs have been
given by Kalman~\cite{DKM}, and Minda and Phelps \cite{MP}. Maxime
B\^{o}cher proved the theorem in 1892 and then
John~H.~Grace proved it in 1902.  Surprisingly, a significant
generalization was proved prior to both B$\hat{\mbox{o}}$cher and
Grace by Siebeck in 1864. He showed that the critical points of an
$n^{th}$ degree polynomial are the foci of the curve of class $n-1$
which touches each line segment, joining the roots of the
polynomial, at its midpoints. The B\^{o}cher-Grace
Theorem is the $n=3$ case of Siebeck's Theorem.  Marden's book
\cite{M} gives a wonderful introduction to this material and it has
an extensive bibliography.

The purpose of this paper is to give a new and different
generalization the B$\hat{\mbox{o}}$cher-Grace Theorem.

\section{Background}

A \textit{conformal similarity transformation} in the complex plane $\mathbb{C}$ is a complex function that
takes the form
\begin{eqnarray}\label{conformalDefinition}
S(z)=\alpha +\beta \,z
\end{eqnarray}
for complex numbers $\alpha$ and $\beta\ne0$.  Upon writing
$\beta=re^{i\theta}$ it is easy to see what action is taken upon $z$
when applying $S$: a rotation, a uniform scaling of the real and
imaginary parts of $z$, and a translation in the plane.  It is
understood that conformal similarity transformations preserve the
eccentricity of ellipses, and for this reason any ellipse in the
complex plane can be mapped via a conformal similarity
transformation to an ellipse with foci $\pm1$.

When studying relationships between the critical points of a
polynomial and its roots, conformal similarity transformations play
a very important standardizing role. Specifically, if $S$ is a
conformal similarity transformation as defined in
equation~(\ref{conformalDefinition}) and a complex polynomial
$p(z)=\prod_{k=1}^n(z-z_k)$ has critical points $r_k$,
$k=1,\ldots,n-1$, then the polynomial
$P(z)=\prod_{k=1}^n\big(z-S(z_k)\big)$ has critical points $S(r_k)$,
$k=1,\ldots,n-1$.  A proof of this fact can be found in \cite{DKM}.

An \textit{affine transformation} in $\mathbb{C}$ is a complex function that takes the form
\begin{eqnarray}\label{affineGeneral}
\Phi(z)=\alpha \,z +\beta \,\overline{z}+\gamma
\end{eqnarray}
for complex numbers $\gamma$, $\alpha\ne0$ and $\beta\ne0$.  Unlike
the conformal similarity transformation it is less clear what action
is taken when applying $\Phi$ to $z$.  However, if
$|\alpha|\ne|\beta|$ then we can introduce real parameters $a>|b|>0$
and $\varphi,\theta\in(-\pi,\pi]$ so that
$$
\alpha=\frac{1}{2}(a+b)e^{i(\varphi+\theta)}\,\mbox{\text{ and }}\,
\beta=\frac{1}{2}(a-b)e^{i(\varphi-\theta)}.
$$
With these new parameters $\Phi$ can be decomposed into three distinct components: a rotation, a
purely affine transformation, and a conformal similarity transformation.  That is,
\begin{eqnarray}\label{affineDecomposed}
\Phi(z)=(S\circ A\circ R)(z)
\quad\mbox{\text{ where }}
\begin{cases}
R(z)= \displaystyle {e^{i\theta}z},\\
A(z)=\displaystyle{\frac{a+b}{2c}\,z+ \frac{a-b}{2c}\,\overline{z}},\\
S(z)= \displaystyle {c\,e^{i\varphi}z+\gamma},
\end{cases}
\end{eqnarray}
where $c=\sqrt{a^2-b^2}$.  We refer to $A$ as a purely affine
transformation because it is the component of $\Phi$ that
distinguishes affine from conformal similarity transformations,
permitting independent scaling of the real and imaginary parts of
$z$.  It is easy to verify, from equation~(\ref{affineDecomposed}),
that
\begin{eqnarray}
A(e^{it})=\frac{a}{c} \cos t+i\frac{b}{c}\sin t,
\end{eqnarray} and so the image of the unit circle, parametrized by $e^{it}$,
under the affine transformation $A$ is an ellipse with eccentricity
$c/a$ and foci $\pm1$. Additionally, we note that the images of the
rotated roots of unity $e^{i(\theta+2k\pi/n)}$, $k=1,\ldots,n$,
under $A$ are
\begin{eqnarray}\label{rotatedRootImages}
\frac{a}{c} \cos \Big(\theta+\frac{2k\pi}{n}\Big)+i\frac{b}{c} \sin \Big(\theta+\frac{2k\pi}{n} \Big),
\quad k=1,\ldots,n.
\end{eqnarray}

We require several properties of affine transformations,
properties we summarize here, without proof.
Affine transformations are invertible if $|\alpha|\ne|\beta|$.
Thus not only is the image of the unit circle under an affine transformation an ellipse,
but ellipses in $\mathbb{C}$ can be mapped onto the unit circle via an affine transformation.
Affine transformations preserve parallel lines, and preserve the midpoints of line segments.


\begin{figure}[h!]
\centerline{\includegraphics[scale=.75]{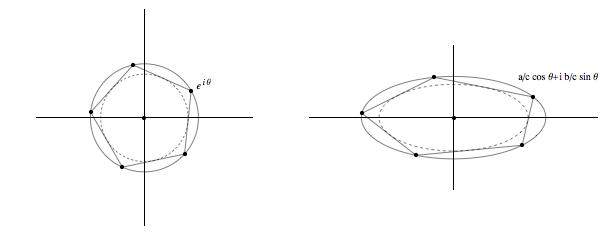}}
\caption{Rotated roots of unity; Their affine image under $A$.}
\label{FigAffinelyRegularPolygon}
\end{figure}

We shall call a convex $n$-gon $\cal{P}$ \textit{affinely regular} if it is the affine image of the regular convex
polygon formed by the $n$ roots of unity.  In Figure \ref{FigAffinelyRegularPolygon}  we illustrate an
affinely regular polygon in the case when $n=5$.  At left is the application of $R(z)=e^{i\theta} z$ to the
roots of unity $e^{i2k\pi/n}$, $k=1,\ldots,n$; at right is the application of $A$ to these rotated roots of unity;
here $S(z)=z$. We also draw attention to the circumscribing and inscribing circles at left, and the corresponding
circumscribing and inscribing ellipses at right.  The circles have radii 1 and $\cos\pi/n$, while the ellipses have
foci $\pm1$ and $\pm \cos\pi/n$, respectively.  In light of the summarized properties of affine transformations,
we note that since the inscribing circle is tangent at the midpoints of the regular polygon pictured, the inscribed
ellipse is also tangent at the midpoints of the corresponding affinely regular polygon.

In fact the existence of an inscribed ellipse interpolating the midpoints of a convex polygon can be shown to
characterize affinely regular polygons.  To see this, let $\cal{P}$ be a polygon with an inscribed ellipse
interpolating its midpoints, and let $\cal{P'}$ denote the convex polygon formed by the midpoints.
Let $\Phi$ be an affine transformation that maps the inscribed ellipse to the unit circle.
The images of the vertices of $\cal{P'}$ under $\Phi$ lie on the unit circle, and the sides of
$\Phi(\cal{P})$ are tangent to the unit circle. There are $2n$ right triangles that can be formed using the origin,
a vertex from $\Phi(P')$, and an adjacent vertex from $\Phi(\cal{P})$.  It is not difficult to see that each of these
triangles is congruent to all the other such triangles.  Hence, the vertices of $\Phi(\cal{P})$ all lie on a common
circle, implying that $\cal{P}$ is an affinely regular polygon.


Here we list some useful properties of the Chebyshev polynomials that can be found in \cite{MH}.
The \textit{Chebyshev polynomials} $T_n(x)$ of the first kind and $U_n(x)$ of the second kind can be defined by
\begin{eqnarray}\label{realCheb}
T_n(\cos\theta)=\cos( n \theta),\mbox{\text{ and }}U_{n}(\cos\theta)=\frac{\sin (n+1) \theta}{\sin \theta },\quad n\ge0.
\end{eqnarray}
From equation~(\ref{realCheb}) we see that $T_n'(x)=nU_{n-1}(x)$,
and that the $n-1$ roots of $U_{n-1}(x)$ are $\{\cos k\pi/n:
k=1,\ldots,n-1\}$.  For complex $z$ the formula for $T_n$ is given
by
\begin{eqnarray} \label{complexDef}
T_n(z)=\frac{1}{2}\bigg((z+\sqrt{z^2-1})^n+(z-\sqrt{z^2-1})^n\bigg).
\end{eqnarray}

From our decomposition of affine transformations in equation
(\ref{affineDecomposed}), a general affine transformation can be
thought of first as a mapping of the unit circle to an ellipse with
foci $\pm1$, followed by a conformal similarity transformation.
There is effectively a one-to-one correspondence between the
functions $A$ and the family of ellipses with foci $\pm1$.  This
family, parametrized by $s>0$, is given by
\begin{eqnarray}\label{family}
\cosh s\cos t+i\sinh s\sin t, \quad 0\le t<2\pi.
\end{eqnarray}
The foci of each ellipse is $\pm1$ since $\cosh^2s-\sinh^2s=1$. The
eccentricity of each ellipse is $0<\sech s<1$.  Since conformal
similarities preserve the eccentricity of an ellipse, it follows
that any ellipse in the plane can be mapped by a conformal
similarity transformation to a member of this family.

A property of the Chebyshev polynomials $T_n(z)$ that deserves to be better known is that they are periodic on
 the
family of ellipses parametrized in equation (\ref{family}). We explain this assertion in the following lemma.

\medskip\noindent {\bf Lemma:} Let $a>b>0$ with $c^2=a^2-b^2$. The complex function
\begin{eqnarray}
f(t)=T_n\bigg(\frac{a}{c}\cos t+i \frac{b}{c}\sin t\bigg)
\end{eqnarray}
is periodic for real $t$, with period $2\pi/n$.  Said another way, the Chebysev polynomial
$T_n(z)$ is periodic on the ellipse
\begin{eqnarray}\label{ellipse}
\frac{a}{c}\cos t+i \frac{b}{c}\sin t,
\end{eqnarray}
taking on the value $\displaystyle{T_n\bigg(\frac{a}{c}\cos \theta+i \frac{b}{c}\sin \theta\bigg)}$
exactly $n$ times.

\medskip
\noindent{\bf Proof:} Key to this result is the observation that if $a^2-b^2=1$ and if
$z=a \cos t +ib \sin t $, then
\begin{eqnarray}\label{connection}
\sqrt{z^2-1}= b \cos t +ia \sin t
\end{eqnarray}
(see exercise 5, page 17 in \cite{MH}).  Equation (\ref{connection}) implies that
\begin{eqnarray}\label{zpm}
z~\pm~\sqrt{z^2-1}=(a~\pm~b)\cos t+i(a~\pm~b)\sin t.
\end{eqnarray}
Substituting equation~(\ref{zpm}) into the formula
(\ref{complexDef}) yields
\begin{eqnarray}
\nonumber T_n(a \cos t+ib \sin t)=\frac{1}{2}\bigg((a + b)^n e^{i n
t} + (a - b)^n e^{-i n t}\bigg).
\end{eqnarray}
Of course, $e^{int}$ is periodic with period $2\pi/n$.

More generally, if $a^2-b^2=c^2\ne1$, then
\begin{align}\label{theLink}
T_n\bigg(\frac{a}{c} \cos
\Big(\theta+\frac{2k\pi}{n}\Big)+i\frac{b}{c} \sin
\Big(\theta+\frac{2k\pi}{n} \Big)\bigg)=\frac{1}{2}
\bigg(\Big(\frac{a + b}{c}\Big)^n e^{i n \theta} + \Big(\frac{a -
b}{c}\Big)^n e^{-i n \theta}\bigg),\end{align} for $k=1,\ldots,n$.
$\blacksquare$

\medskip
At this stage we invite the reader to compare the formula for the
affine image of the rotated roots of unity under the affine
transformation $A$ in equation~(\ref{rotatedRootImages}) with the
argument of $T_n$ in equation~(\ref{theLink}).  They are the same.
This means that $T_n$ is constant on the affine images of the
rotated roots of unity.  This is true for every affine
transformation $A$, as defined in equation~(\ref{affineDecomposed}).

\begin{figure}[h!]
\centerline{\includegraphics[scale=.75]{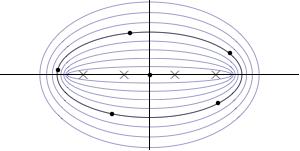}}
\caption{Ellipses with foci $\pm1$; Affine image of rotated roots of
unity.}\label{ellipticalFamily}
\end{figure}

In Figure \ref{ellipticalFamily} we illustrate some members of the family (\ref{family}).
In addition, we plot the $n=5$ images of a particular rotation of the $n$ roots of unity under
a particular affine transformation $A$.  The polynomial $T_n$ is constant on these $n$ points.
Along the real axis we also mark the $n-1$ roots of $T'_n(z)=nU_{n-1}(z)$

\section{Main Result}
\medskip
\noindent{\bf B$\hat{\mbox{\text{o}}}$cher-Grace Theorem for
Polygons:} Let $p$ be an $n^{th}$ degree complex polynomial and let
its critical points take the form
\begin{eqnarray}\label{criticalPoints}
\alpha+\beta \cos k\pi/n, \quad k=1,\ldots,n-1, \quad\beta\ne0.
\end{eqnarray}
There is an inscribed ellipse interpolating the midpoints of the convex polygon formed by the roots of
$p$, and the foci of this ellipse are the two most extreme critical points of $p$: ~$\alpha\pm\beta \cos \pi/n.$

\medskip
\noindent{\bf Proof:}  Assume that the critical points of $p$ take
the form in expression~(\ref{criticalPoints}). Without loss of
generalization, exploiting properties of conformal similarity
transformations, we may assume $\alpha=0$ and $\beta=1$. Thus
$p'(z)=\gamma U_{n-1}(z)$, and $p(z)= \gamma/n T_n(z)+\delta$, for
complex constants $\gamma\ne0 $ and $\delta $.  If we designate  by
$z_k$, $k=1,\ldots,n$, the roots of $p$, then
$T_n(z_k)=-n\delta/\gamma $, $k=1,\ldots,n$.  That is to say, $T_n$
is constant on the $z_k$, $k=1,\ldots,n$.   As noted in our
discussion, $T_n(z)$ can only be constant for $n$ distinct points in
$\mathbb{C}$, all of which lie on a single ellipse from the family
(\ref{family}).  Thus the $z_k$ must take the form
(\ref{rotatedRootImages}). In turn this implies that the roots of
$p$ form an affinely regular polygon $\cal{P}$, and hence admit an
inscribed ellipse interpolating the midpoints of $\cal{P}$. In our
discussion of affinely regular polygons, we noted that the foci of
the inscribed ellipse are $\pm \cos\pi/n$, consistent with the
assertion of the theorem. $\blacksquare$

\medskip
When $n=3$ the critical points \textit{always} take the form $\alpha\pm \beta\cos \pi/3$.
That is, all triangles are affinely regular; all triangles can be mapped to an equilateral
triangle with an affine transformation. This is the original B$\hat{\mbox{\text{o}}}$cher-Grace theorem.
When $n=4$ the critical points take the form $\alpha\pm \beta\cos k\pi/4$, $k=1,2,3$, if and only if the roots of
$p$ form a parallelogram.  Said another way, the only quadrilaterals that are affine images of squares are
parallelograms.

The form of the critical points in expression~(\ref{criticalPoints})
can be used to characterize affinely regular polygons. That is, a
polygon $\cal{P}$ formed from the roots of a polynomial $p$ is
affinely regular if and only if the critical points of $p$ take the
form expression~(\ref{criticalPoints}).  With this observation we
conclude this article with a stronger statement than the generalized
B$\hat{\mbox{\text{o}}}$cher-Grace Theorem above.

\medskip
\noindent{\bf B$\hat{\mbox{\text{o}}}$cher-Grace Characterization
Theorem:} Let $p$ be an $n^{th}$ degree polynomial and let $\cal{P}$
denote the convex polygon formed by the roots of $p$.  The polygon
$\cal{P}$ admits an inscribed ellipse interpolating its midpoints if
and only if the critical points of $p$ take the form of the
expression expression~(\ref{criticalPoints}).

\end{spacing}

\end{document}